\documentclass[12pt,a4paper]{article}

    \usepackage{amsmath, amsthm, amssymb}

 \topmargin= - 1 true mm
\renewcommand{\b}[1]{\mbox{\boldmath $#1$}}

\newcommand{\sgn}{\operatorname{{\mathrm sgn}}}
\newcommand{\sign}{\operatorname{{\mathrm sign}}}

\newcommand{\exn}{\mathbb{E}}
\newcommand{\pr}{\mathbb{P}}
\newcommand{\Nn}{\mathbb{N}}
\newcommand{\Rr}{\mathbb{R}}
\newcommand{\Ac}{\mathcal{A}}
\newcommand{\ep}{\varepsilon}
\newtheorem{theorem}{Theorem}
\newtheorem{corollary}[theorem]{Corollary}

\newtheorem{example}[theorem]{Example}

\begin{document}

\title{On the Expectations of Maxima of Sets of Independent Random Variables}

\author{D.\ V.\ Tokarev\footnote{ARC Centre of Excellence for Mathematics and Statistics of Complex Systems
(MASCOS). E-mail: {\tt daniel.tokarev@sci.monash.edu.au} (the corresponding author).}
 \ %
 and K.\ A.\ Borovkov\footnote{Department of Mathematics and Statistics, University of
Melbourne, Parkville 3010, Australia.}}

\date{}
\maketitle

\begin{abstract}
Let $X^1, \ldots, X^k$ and $Y^1, \ldots, Y^m$ be jointly independent copies of random
variables $X$ and $Y$, respectively. For a fixed total number $n$ of random variables,
we aim at maximising $M(k,m):=\mathbb{E} \max \{X^1, \ldots, X^k, Y^1, \ldots, Y^{m}
\}$ in $k = n-m\ge 0$, which corresponds to maximising the expected lifetime of an
$n$-component parallel system whose components can be chosen from two different types.
We show that the lattice $\{M(k,m)$: \ $k, m\geq 0\}$ is concave, give sufficient
conditions on $X$ and $Y$ for $M(n,0)$ to be always or ultimately maximal and derive a
bound on the number of sign changes in the sequence $M(n,0)-M(0,n)$, $n\geq 1$. The
results are applied to a mixed population of Bienayme-Galton-Watson processes, with the
objective to derive the optimal initial composition to maximise the expected time to
extinction.
\end{abstract}

{\noindent\small{\em Key words and phrases.} Parallel systems, expected lifetime, branching processes.\\
 \vspace{5mm}%
}

\section{Introduction}

Consider a simple $n$-component parallel system (called an \textit{$n$-assembly} in
what follows), whose component lifetimes are independent non-negative random variables
$X_1, \ldots, X_n$ with finite means. The $n$-assembly fails when all the $n$
components fail, and so its expected lifetime is given by $\exn \max\{X_1,X_2,\ldots,
X_n\},$ which is a standard measure of the system's reliability.

Such quantities and related characteristics of complex systems often appear in
literature on reliability. For instance, comparisons and bounds for expected lifetimes
of series and parallel systems are considered in \cite{MP}, while assemblies whose
components' lifetime distributions can depend on the amount of ``investment" are
discussed in~\cite{JA}. Expected lifetimes conditional on survival up until time $t$
are investigated in~\cite{AB}, while assemblies of dependent components with
interchangeable lifetime distributions are analysed in~\cite{NRS}. More recently, upper
and lower bounds for the expected lifetimes of $n$-assemblies of components with
(possibly) non-identically distributed lifetimes in terms of the expected lifetimes of
assemblies with identically distributed components' lifetimes were given
in~\cite{HJST}.

In the present study, we assume that, when building an $n$-assembly, one may choose
components from $d$ different types that have different lifetime distributions. Should
one choose all $n$ components of the same type ({\em unmixed\/} assembly) or should one
mix? In general, the answer will depend not only on the distributions of the random
lifetimes of the components, but also on $n$. In this paper we shall examine these
relationships.

Questions of such kind arise naturally in reliability problems such as server failures,
where information is duplicated on multiple hard disks. In conservation biology, the
``Single Large or Several Small'' problem provides another application. Here one should
decide between the creation of a single large protected habitat or several small ones.
This may be modelled by $n$-assemblies since the lifetimes of individuals in the
species in different small habitats may be modelled by random variables $X_i$ with
distinct (due to differing environmental conditions) distributions. In contrast, a
single large habitat with the same lifetime distribution for all individuals in the
species corresponds to an $n$-assembly composed of identically distributed random
variables. Furthermore, the question of whether mutation is inherently advantageous
because it leads to diversity may also be analysed in a similar way.

In this paper, we first show that the expected lifetime $M(k_1, \ldots, k_d)$ of an
assembly consisting of $k_i$ components of type $i=1, \ldots, d$, is a concave function
of $(k_1, \ldots, k_d)$. This implies that one can easily find the optimal $n$-assembly
in the general case using the discrete steepest ascent method in linear time in $dn$.
After that we focus on the case of two component types, with a view to maximising the
life expectancy of the system, noting that generalisations to more than two types are
straightforward and are omitted for the sake of simplicity. We denote the two random
lifetimes by $X$ and $Y$ and their respective distribution functions by $F$ and $G$,
always assuming that $\exn X+\exn Y<\infty$, and assume that $X_1, X_2, \ldots$ are
i.i.d.\ like $X$,  and $Y_1, Y_2,\ldots$ are i.i.d.\ like $Y$, with all the $X_i$'s and
$Y_j$'s being jointly independent. Our main object of study is
\[
M(k, m):=\exn \max \{X_1,\ldots, X_k,Y_1,\ldots, Y_m\}
 =\int_0^\infty \bigl(1- F(s)^k G(s)^{n-k}\bigr) ds.
\]
We classify the possible behaviour of these quantities, e.g.\ when it is advantageous
to always choose one type and when it is the right choice eventually, giving sufficient
conditions for these. We also give a simple bound for the number of sign changes of
$M(n,0)-M(0,n)$ and apply the aforementioned results to the case of
Bienayme-Galton-Watson processes, with a view to selecting an optimal initial
composition in order to maximise the expected time to extinction of the entire
population. But first we present a couple of simple examples that show some possible
types of behaviour that can be displayed by $M(k,m)$.

\begin{example}
{\rm Suppose we can choose from two types of components: one with lifetimes $X\sim
U[0,1]$ (i.e.\ distributed uniformly on $[0,1]$) and   the other with deterministic
lifetimes $Y\equiv \ep =\text{const}$. Clearly, if $\ep \geq 1$ then one should always
choose components of type~$Y$. In general, we will refer to a situation where the
optimal (with the longest expected lifetime) $n$-assembly is always composed of one
particular type of component as the  {\em dominant\/} case, and to that component type
as the {\em dominant\/} type.

Now assume that $0<\ep <1$. The choice is between selecting all $n$ components of the
first type or $n-1$ components of the first type and one component of the second type
(there is no point in taking more than one type two component to extend the expected
lifetime as all such components will fail at precisely the same time). Integrating to
obtain the corresponding expectations gives
\begin{align*}
M(n,0)
 &=\int_0^{1}\big(1-x^n\big)\, dx=1-\frac{1}{n+1},
  \\
M(n-1,1)
 &=\ep + \int_{\ep }^1 \big(1-x^{n-1}\big)\,dx=1-\frac{1}{n}+\frac{\ep ^n}{n}.
\end{align*}
We conclude that the assembly corresponding to $(n,0)$ yields a higher expected
lifetime than that for $(n-1,1)$ for $n$ such that $(n+1)^{1/n}<\ep ^{-1}$. If $\ep
<1/2$, this will always be true, and then we have a dominant case. Since $\ep <1$ and
$(n+1)^{1/n}\to  1$ as $n \to  \infty$, we see that, regardless of $\ep $,
$M(n,0)>M(n-1,1)$ holds for all large enough $n$. In general, if unmixed assemblies
$\{X_1, \ldots, X_n\}$ are advantageous for all large enough $n$, we say that $X$ is
{\em ultimately dominant}. We will see (in Theorem~\ref{t6}) that this behaviour takes
place whenever one distribution function ultimately dominates the other.

Since $(n+1)^{1/n}$ is a  decreasing function of $n\geq 1$, we see that the values of
$n$ for which $(n+1)^{1/n}>\ep ^{-1}$, i.e.\ for which $(n-1,1)$ is preferable to
$(n,0)$, form either an empty set or a finite sequence of successive integers. In the
latter situation, since the optimal $n$-assemblies are composed of more than one
component type, we say that we have a  {\em non-dominant case\/} with the number of
{\em dominance changes\/} (i.e.\ sign changes of $M(n,0)-M(0,n)$) equal to $1$. }
\end{example}

We will see that this number admits a simple upper bound in terms of the distribution
functions of $X$ and $Y$ (Theorem~\ref{t7}), but can actually be infinite, as the
following example suggested by A.~Sudbury shows.

\begin{example}\label{e2}
{\rm Let $X$ and $Y$ have distributions $\pr (X=2i)=p_i$ and $\pr (Y=2i+1)=q_i,$
$i=0,1,2, \ldots, $ such that
\[
P_j:=p_j+p_{j+1}+\cdots=q_{j+1}+q_{j+2}+\cdots=2^{-2^j}, \qquad j\geq 2.
\]
Set $\Ac (k,m):=\max\{X_1,\ldots,X_k,Y_1,\ldots, Y_m\}$, so that $M(k,m)=\exn \Ac
(k,m).$

We will show that there exist  two increasing infinite sequences $\{k_i\}$ and
$\{l_i\}$ such that $M(k_i,0)>M(0,k_i)$ and $M(l_i,0)<M(0,l_i)$ for all sufficiently
large $i$.

Indeed, taking $k_i:=\lceil 2^{2^{i}} i^{-2}\rceil$, where $\lceil x \rceil =
\min\{k\in \mathbb{Z}: \, k\geq x\}$, we see that
\begin{align*}
&\pr \bigl(\Ac (k_i,0)<2i\bigr) =(1-P_i)^{k_i}=(1-2^{-2^i})^{2^{2^i}  i^{-2}}
 = \exp\bigl\{-(1+o(1))i^{-2}\bigr \}\to 1,
  \\
&\pr \bigl(\Ac (k_i,0)<2(i-1)\bigr)=(1-2^{-2^{i-1}})^{2^{2^i} i^{-2}}
 = \exp\bigl\{-(1+o(1))2^{2^{i-1}}i^{-2}\bigr\}\to 0
\end{align*}
as $i \to  \infty$, so that $M(k_i,0)\sim 2(i-1)$ (here $\sim$ means asymptotic
equivalence: $u_n\sim v_n$ as $n\to  \infty$ iff $u_n/v_n \to  1$), whereas
\begin{align*}
 \pr \bigl(\Ac (0,k_i)\geq 2i+1\bigr)
  &=1-(1-P_{i-1})^{k_i}\\
  &=1-(1-2^{-2^{i-1}})^{2^{2^i} i^{-2}}=1-\exp\{-(1+o(i))2^{2^{i-1}}i^{-2}\} \to 1.
\end{align*}
Thus $M(k_i,0)<M(0,k_i)$ for all sufficiently large $i$.

Now let $l_i:= 2^{2^{i+1}} (i+1)^{-2}$. Then it can be shown in a similar fashion that
$M(l_i,0)>M(0,l_i)$ for all sufficiently large $i$, so that we have infinitely many
dominance changes in this example.}
\end{example}

\section{Concavity of the expected lifetimes and\\ dominance}

When  maximizing the function $M(k_1, \ldots, k_d),$ it is helpful to know that it is
concave. In particular, this guarantees that, moving in the direction of the steepest
ascent one can find the optimal composition in linear time in~$dn$.

\begin{theorem}
 \label{t3}
Let $X_1, \ldots, X_d$ be random lifetimes with corresponding distribution functions
$F_1, \ldots, F_d$. Then the function
\[
M (\b{x}) : =\int_0^\infty \big(1-F_1(s)^{x_1}\cdots F_d(s)^{x_d}\big)\, ds, \qquad
 \b{x}=(x_1,\ldots, x_d)\in \Rr_+^d,
\]
is concave.
\end{theorem}

Clearly, the function $M (\b{x})$ coincides with the corresponding systems' expected
lifetimes at the integer lattice points.

\begin{proof} We will begin with the obvious observation that, for any $\b{x}\in \Rr_+^d$,
$$
 M (\b{x})\leq M (\lceil x_i\rceil,\ldots,\lceil x_d\rceil)\leq \sum_{j=1}^d \lceil x_j\rceil \exn  X_j<\infty.
$$
Fix arbitrary $\b{x},\ \b{y}\in \Rr_+^d$ and let $F(s):=F_1^{x_1}(s)\cdots
F_d^{x_d}(s),$ $G(s):=F_1^{y_1}(s)\cdots F_d^{y_d}(s)$. It suffices to show that the
univariate function
$$
M _1(z):=M (z\b{x}+(1-z)\b{y})=\int_0^\infty \big(1-F(s)^zG(s)^{1-z}\big)\, ds,
 \qquad z\in [0,1],
$$
is concave. To justify differentiation inside the integral, we first fix
an arbitrarily small $\delta>0$ and show that the integral of the integrand's
derivative converges absolutely and uniformly in $z\in(\delta,1-\delta)$. We have
$$
\int_0^\infty \Big|\frac{\partial}{\partial z}\big(1-F(s)^z G(s)^{1-z}\big)\Big|\,ds
 =\int_0^\infty\big|\ln F(s)-\ln G(s) F(s)^zG(s)^{1-z}\big|\, ds,
$$
where the integrand is dominated by the function
$$
|\ln F(s)|F(s)^\delta + |\ln G(s)|G(s)^\delta
$$
(using the natural convention  $0\cdot \ln 0=0$), which is clearly bounder on
$(0,\infty)$ and so is integrable there since $|\ln F (s)| F(s)^\delta\sim 1-F(s)$ as
$s\to\infty$.

In a similar way one can  verify that $M _1$ is twice differentiable and
\begin{align*}
M _1''(z)
 &= \int_0^\infty\frac{\partial^2}{\partial z^2}\big(1-F(s)^zG(s)^{1-z}(s)\big)\, ds\\
 &=-\int_0^\infty \big(\ln F(s)-\ln G(s)\big)^2F(s)^zG(s)^{1-z}\, ds<0
\end{align*}
unless $F(s)\equiv G(s)$. This shows that $M_1$ and hence $M$ are concave.
\end{proof}

Now we turn our attention to the case when  one type is always dominant, i.e.\ one has
$M(n,0)\geq M(k,n-k)$  for any $n\in \Nn,$  $0\leq k \leq n$.

It is obvious that if
\begin{equation}
 \label{SD}
 F(s)\leq G(s)  \qquad  \text{ for all\ } s\in \Rr_+,
\end{equation}
then we have a dominant case. One can weaken $\eqref{SD}$ by requiring that one
distribution always has a heavier integrated tail than the other.

\begin{theorem}{\rm(Theorem~\ref{t7}.6 of \cite{BP})}
 \label{t4}
Let $X$ and $Y$ be random lifetimes with df's $F$  and~$G$ respectively. If
\begin{equation}
 \label{ICO}
 \int_{s}^\infty \left(F(t)-G(t)\right) dt<0
  \qquad \text{for all \ } s\in \Rr_+,
\end{equation}
then $M(n,0)>M(k,n-k)$ for all $n \in \Nn$ and $k\in[ 0, n)$.
\end{theorem}

Recall that \eqref{ICO} defines the so-called increasing convex ordering on the set of
distributions on $\Rr_+$ (see e.g.\ p.11 in~\cite{Sz}), whereas  \eqref{SD} is the
standard stochastic ordering.

The next example shows that the above condition can be necessary in some special cases.

\begin{example}
 \label{e5}
{\rm Let $X$ be uniformly distributed on $[0,1]$ and $Y$ uniformly distributed on
$[a,a+\ep ]$ with $0<a<a+\ep <1$. Then clearly
\begin{align*}
 M(n,0)-M(0,n)
  &=\int_0^1(1-x^n)\, dx-\left[a+\int_a^{a+\ep }\bigg(1-\Big(\frac{x-a}{\ep }\Big)^n\bigg)\, dx\right]\\
  &=(1-\ep )\bigg(1-\frac{1}{n+1}\bigg)-a.
\end{align*}
Hence $M(n,0) \geq M(0,n)$ iff ${a}/{(1-\ep)}\leq 1- {(n+1)^{-1}}$. This is the case
for all $n$ (and then the first type is dominant) iff ${a}/{(1-\ep)} \leq  {1}/{2},$
and so this condition is necessary for the first type to dominate. On the other hand,
by Theorem~\ref{t4} for the $X$ to be dominant it suffices that~\eqref{ICO} holds,
which is equivalent to the condition that the area of the upper of the two triangles
formed by the graphs of the functions $F$ and $G$ is no less than the area of the lower
one. Since the two triangles are similar, that means that the abscissa of the point of
the intersection of the two graphs does not exceed $1/2$, i.e.\ ${a}/{(1-\ep)}\leq
{1}/{2}$. Thus we see that condition~\eqref{ICO} which is sufficient for dominance of
$X$ is equivalent to the condition that was shown to be necessary for the dominance.}
\end{example}

There is also a simple sufficient condition which ensures that one type will eventually be dominant.

\begin{theorem}
 \label{t6}
Let $X$ and $Y$ be random lifetimes with the corresponding distribution functions $F$
and $G$, and let there exist an $s_0 \in \Rr_+$ such that  $F(s)> G(s)$ for all $s\geq
s_0$. Then there exists an $n_0\in \Nn$ such that, for all $n>n_0$ and $1\leq k\leq n$,
$M(k,n-k)< M(0,n)$, i.e.\ $Y$ is ultimately dominant.
\end{theorem}

\begin{proof}
We will prove the result in the continuous case, an extension to the general case being
obvious.

Let $s_0=\inf\{s: F(t)>G(t), \ t>s\}$ and $F(s_0)=G(s_0)=:e^{-a}$. Then
\begin{align*}
M(k,n-k)-M(n,0) & =
 \int_0^\infty \big(F(s)^kG(s)^{n-k}-G(s)^n\big)\, ds\\
 &= \int_0^\infty G(s)^n\left(\left(\frac{F(s)}{G(s)}\right)^k-1\right)\, ds\\
 &= \int_0^{s_0}+\int_{s_0}^\infty=:I_1 + I_2.
\end{align*}
Clearly,   $|I_1  |\leq s_0e^{-an}$. Now set $H(s):=\ln {F(s)}-\ln {G(s)}\geq 0$ for
$s>s_0$ and choose $s_1>s_0$ such that $G(s_1)=:e^{-b}>e^{-a}$. Then
\[
I_2 =\int_{s_0}^\infty G(s)^n\left(e^{kH(s)}-1\right)\, ds \geq \int_{s_1}^\infty
G(s)^n\left(e^{kH(s)}-1\right)\geq ke^{-bn}\int_{s_1}^\infty H(s)\, ds.
\]
Since  $a - b>0$, one has  $\displaystyle{\int_{s_1}^\infty H(s)\, ds >
s_0e^{-(a-b)n}}$ for all large enough~$n$, which clearly implies that, for all $k\geq
1,\ n\geq n_0,$
$$
I_2\geq ke^{-bn}\int_{s_1}^\infty H(s)\, ds>s_0e^{-an}\geq |I_1|,
$$
as required.\end{proof}

Using Theorem~\ref{t6}, one can easily find examples of non-dominant cases. It suffices
to have $\exn X>\exn Y$ and $F(s)>G(s)$ for $s>s_0\in \Rr_+$, so that
$M(1,0)-M(0,1)>0,$ but $M(n^*,0)-M(0,n^*)<0$ for some $n^*\in \Nn$ (in other words, at
least one sign change of $M(n,0)-M(0,n)$ occur between $n=1$ and $n=n^*$). In such
situations, there typically exist $n$ and $k^*\in (0,n)$ such that $M(k^*,n-k^*)\geq
M(k,n-k)$ for all $0\leq k\leq n$, i.e.\ the mixed assembly $(k^*,n-k^*)$  is optimal.

\section{The number of sign changes of $M(n,0)-M(0,n)$}

One can obtain a useful upper bound for the number of sign changes in the sequence
$M(n,0)-M(0,n)$, $n\geq 1$, by generalizing the Descartes rule of signs. Clearly, this
number does not exceed the number of zeros of the function $$\psi(x):=\int_0^\infty
\big(1-F(s)^x\big)\, ds - \int_0^\infty \big(1-G(s)^x\big)\, ds, \qquad x>0.$$ In order
to find a bound for this number, we will employ a generalisation of the Descartes rule
of signs~\cite{CP}.

One says that an ordered system $\{u_1,u_2,\ldots\}$ of functions defined on
$I\subseteq \Rr$ satisfies the Descartes rule of signs on $I$ if the number of zeros
(with multiplicities) of a linear combination $\sum a_ju_j(x)$ of these functions is
less than or equal to the number of variations of strict sign in the sequence $\{a_1,
a_2, \ldots\}$ of the coefficients in the combination.

A system of functions $\{u_1, \ldots, u_n\}$ defined on $I\subseteq \Rr$ is said to be
strictly sign regular if, for any fixed $k\leq n$ and any $x_1<x_2<\cdots<x_m$, $m\geq
k$, with all $x_j\in I, \ j\leq m$, all minors of order $k$ of the matrix
$\big(u_j(x_i)\big)_{i\leq m,\, j\leq n}$ have the same strict sign. For example,
systems of the form $\{x^{\beta_1}, \ldots, x^{\beta_n}\}$ and $\{e^{-\beta_1x},\ldots,
e^{-\beta_nx}\}$ with $0<\beta_1<\beta_2<\cdots<\beta_n$ are strictly sign regular on
$\Rr$ (see e.g.\ p.34 of~\cite{CP}). Proposition 2.6 of \cite{CP} asserts that a system
of functions $\{u_1, \ldots, u_n\}$ satisfies the Descartes rule of signs iff it is
strictly sign regular. \vspace{0.3cm}

Given two distribution functions $F$ and $G$ on $\Rr_+$,  let
$$
 \mu(A):=\int_A \big[F^{-1}(e^{-t})-G^{-1}(e^{-t})\big]\, dt
$$
be a signed measure on $\Rr_+$. Let $S^+$ be a positive set and $S^-=\Rr_+\setminus
S^+$ a negative set from the Hahn-Jordan decomposition $\mu=\mu^+-\mu^-$. We assume
(for the moment) that $S^+$ has only finitely many connected components $S_1^+, \ldots,
S^+_{\nu_+}$  and that, without loss of generality, $0\in S^+_1$. By $S^-_1, \ldots,
S^-_{\nu_-}$ we similarly denote the connected components of $S^-$, so that
$x_1^+<x_1^-<x_2^+<x_2^-<\cdots$ for any $x_i^\pm \in S_i^\pm$ (note that $0\leq
\nu_+-\nu_-\leq 1$). Let  $\nu=\nu^- +\nu^+$ and introduce sets $Q_1, \ldots, Q_\nu$ by
putting $Q_1:=S_1^+,$ $Q_2:=S_1^-,$ $Q_3:=S_2^+,$ $Q_4:=S_2^- $ etc. The quantity $\nu$
may be called the number of intersections of $F$ and $G$ since $\nu^+$ and $\nu^-$ are
none other than the numbers of connected components of $\{s:F(s)\leq G(s)\}$ and
$\{s:F(s)\geq G(s)\}$, respectively. When $S^+$ has infinite number of connected
components, we set $\nu=\infty$. This was precisely the case in Example~\ref{e2} from
the Introduction.

\begin{theorem}
 \label{t7}
Let $X$ and $Y$ be random lifetimes, $F$ and $G$ being their corresponding distribution
functions. Then the number of sign changes of\/ $M(n,0)-M(0,n)$, $n\geq 1$, does not
exceed $\nu-1$.
\end{theorem}

\begin{proof}
Noting that $X\stackrel{d}{=}F^{-1}(U)$ and $Y\stackrel{d}{=}G^{-1}(U)$, where $U\sim
U[0,1]$ and $H^{-1}(u):=\inf\{t:H(t)>u\}$ denotes the generalized inverse of the
distribution function $H$, we can write
$$
\psi(x)=\int_0^1 F^{-1}(s^{1/x})\, ds - \int_0^1 G^{-1}(s^{1/x})\, ds.
$$
Letting $s^{1/x}=e^{-t},$ we see that $\psi(s)=x\psi_1(x)$, where
\begin{align*}
\psi_1(x)
 &=\int_0^\infty e^{-xt}\,\mu(dt)=\sum_{j=1}^\nu\int_{Q_j}e^{-xt} \, \mu(dt)\\
 &=\sum_{j=1}^\nu (-1)^{j+1}\int_{Q_j}e^{-xt}\, |\mu (dt)|=\sum_{j=1}^\nu a_j\exn e^{-xB_j}
\end{align*}
with $a_j=\mu(Q_j)$ (note that $\sign (a_j)=(-1)^{j+1}$) and $B_j=H_j^{-1}(U_j)$, where
$$
H_j(t) :=\frac{1}{a_j}\int_0^t \mathbb{I}_{Q_j}(s)\, \mu(ds)
$$
are distribution functions supported by the corresponding sets $Q_j$, $j=1, \ldots,
\nu,$ and $\mathbb{I}_A$ denotes the indicator of $A,$ while $U_j\sim U[0,1]$ are
independent random variables.

To show that $\{u_j(x)\}:=\{\exn e^{-xB_j}, j=1,\dots,  \nu\}$ is a strictly sign
regular system of functions on $\Rr_+$, first note that $\{v_j(x)\}:=\{e^{-\beta_jx}\}$
is such a system provided that $0<\beta_1<\beta_2<\cdots<\beta_\nu$  (see p.~34 and
Corollary~3.9 of~\cite{CP}). From this it follows that, for a $k\leq \nu$, the sign of
$\det \big(e^{-\beta_j x_i}\big)_{i, j\leq k}$ is one and the same for any choice of
$0<\beta_1<\cdots <\beta_k$ and $0<x_1<\cdots<x_k$.

Now observe that a $k$th order minor of $\big(u_j(x_i)\big)_{j\leq \nu,\, i\leq m}$
admits an integral representation in terms of the corresponding minors of
$\big(v_j(x_i)\big)_{j\leq \nu,\, i\leq m}$. Indeed, consider without loss of
generality the determinant of $\big(u_j(x_i)\big)_{i, j \leq k}$ for
$0<x_1<\cdots<x_k$. It has the form
\begin{align*}
\sum_\sigma \sgn (\sigma) \prod_{j=1}^k u_j(x_{\sigma(j)})
 &=\sum_\sigma \sgn (\sigma) \prod_{j=1}^k\exn e^{-x_{\sigma(j)}B_j}
  =\exn \sum_\sigma \sgn  (\sigma) \prod_{j=1}^k e^{-x_{\sigma(j)}B_j}\\
 &=\int \cdots \int \sum_\sigma \sgn (\sigma) \bigg(\prod_{j=1}^k e^{-\beta_j x_{\sigma(j)}}\bigg)\,
  dH_1(\beta_1)\cdots dH_k(\beta_k),
\end{align*}
where the summation is over all permutations $\sigma$ of $\{1,\ldots, k\}$,
$\sgn(\sigma)$ is the signature of $\sigma$ and we used the independence of $B_j$. To
complete our argument, it remains to make use of the above remark concerning the sign
invariance of the determinants of the matrix $\big(e^{-\beta_jx_i}\big)$ and recall
that, by construction, we have $\beta_1<\beta_2<\cdots<\beta_k$ in the integrand in the
last displayed formula, due to $\beta_j\in Q_j,$ $j \leq k$.

Thus we have shown that our $\{u_j(x)\}$ is a strictly sign regular system and so, by
the Descartes rule, $\psi_1$ (and therefore $\psi$) can have at most $\nu-1$ zeros in
$(0,\infty)$. Thus $M(n,0)-M(0,n)$ can have at most $\nu-1$ sign changes, as required.
\end{proof}

Note that it follows from Theorem~\ref{t6} that if the (possibly infinite) sequence of
points of intersection of the graphs of the distribution functions $F$ and $G$ is
bounded, the number of dominance changes will necessarily be finite. However, if this
sequence is unbounded (so that $\nu=\infty$), the number of dominance changes can be
infinite as  Example~\ref{e2} in fact shows.

\section{Mean time to extinction of mixed branching populations}

In this section we consider an application of the above results to the optimal
selection for populations of simple branching processes, where $n$-assemblies are
composed of the times to extinction of the progenies of the individuals comprising the
initial population.

Suppose we wish to establish a colony populated by two species, $S_1$ and $S_2$, with
the aim of maximising the colony's expected lifetime. Given a fixed initial total
number $n$ of individuals in the colony, we can choose $k$ individuals of species $S_1$
and $n-k$ individuals of species $S_2$. Next we assume that time is discrete and that
the reproduction of individuals in the colony is governed by independent subcritical
Bienayme-Galton-Watson (BGW) processes specific to the species.

Let $\xi\geq 0$ be an offspring random variable with a distribution $p_i=\pr (\xi=i),$
$i\geq 0,$ and $\xi_{n,i}$  with $n, i=1,2,\ldots$  be independent copies of $\xi$.
Further, for $r\in\Nn$, let $Z^{(r)}(0):=r$  and
$$
Z^{(r)}(n):=\sum_{i=1}^{Z^{(r)}(n-1)}\xi_{n,i},\qquad n\geq 1.
$$
Then $\{Z^{(r)}(n), \, n\ge 0\} $ is a BGW process starting with $r$ individuals (for
$r=1,$ we will simply write $\{Z(n)\}$).

Denote by $T^{(r)}:=\inf\{n\geq 1: \, Z^{(r)}(n)=0\}$   the time to extinction of the
process $\{Z^{(r)}(n)\}$ (we write $T$ for $T^{(1)}$). Observe that
$T^{(r)}=\max\{T^1,\ldots, T^r\},$ where the $T^j$'s are independent copies of $T$. Let
$f(\theta):=\sum_{i=0}^\infty p_i\theta^i$ be the probability generating function (pgf)
of $\xi$ and $f_n(\theta)$ the $n$th functional iterate of~$f$. Then
\begin{equation}
 \label{PT}
\pr (T\leq n)=\pr (Z(n)=0)= f_n(0), \qquad n\ge 0,
\end{equation}
and hence  $\pr (T^{(r)}\leq n)=f_n(0)^r$ and $\exn
T^{(r)}=\sum_{n=0}^\infty\big(1-f_n(0)^r\big)$, where we put $f_0(0)=0$.

Now denote by $\{Z_i(n)\}$ a BGW process describing the growth of a population of type
$S_i$ individuals  (starting with a single individual at time $0$), $i=1, 2,$ and by
$T_i$ the respective times to extinction. We will keep the notation $f, f_n$ for the
offspring pgf and its iterates for species $S_1$ and denote the respective pgf's for
$S_2$ by $g$ and $g_n$. Then the expected time to the extinction of our colony
initially consisting of $k$ individuals of species $S_1$ and $n-k$ individuals of
species $S_2$ is given by
\begin{align*}
M(k,n-k)
  &=\exn \max\{T_1^{(k)},T_2^{(n-k)}\}\\
  &=\exn \max\{T_1^1,\ldots, T_1^k, T_2^1, \ldots, T_2^{n-k}\}
   =\sum_{i=0}^\infty \big(1-f_i(0)^kg_i(0)^{n-k}\big),
\end{align*}
where the $T_i^j$ are independent copies of $T_i$, $i=1,2$. The problem we want to
consider is how to choose $k$ that maximises $M(k,n-k)$ for a fixed $n$. First we
obtain some preliminary facts.

\begin{theorem}
 \label{t8}
For any two subcritical BGW processes $\{Z_1(n)\}$ and $\{Z_2(n)\}$ with
offspring means $\mu_2<\mu_1<1$, there exists an $n_0$ such that $f_n(0)<g_n(0)$ for
all $n>n_0$ and therefore $T_1$ is ultimately dominant, i.e.\ there exists an $n_1$
such that $M(n,0)>M(k,n-k)$ for all $n>n_1$ and $1\leq k\leq n$.
\end{theorem}

\begin{proof}
From \cite{Se} we know that, as $n \to  \infty$,
$$
1-f_n(0)\sim \mu_1^nL_1(\mu_1^n)
 \text{\ \ \ and\ \ \ }
1-g_n(0)\sim \mu_2^nL_2(\mu_2^n),
$$
where $L_i(s)$ are functions slowly varying at $0$, i.e., for any $c>0$, $L_i(cs)\sim
L_i(s)$ as $s \searrow 0$, $i=1, 2$. As is well-known (see e.g.\ Proposition~1.3.6(v)
in~\cite{BGT}), if $L(s)$ is slowly varying at $0$ then, for any $\ep >0$, one has
$s^\ep <L(s)<s^{-\ep }$ for all small enough $s>0$. Therefore, for an arbitrarily small
$\ep $,
$$
\frac{L_1(s_1)}{L_2(s_2)}>(s_1s_2)^\ep
$$
once $s_1$ and $s_2$ are small enough. So, for $0<\ep  <
(\ln\mu_2-\ln\mu_1)/(\ln\mu_2+\ln \mu_1)$ we obtain
$$
\frac{1-f_n(0)}{1-g_n(0)}
 \sim \frac{\mu^n_1}{\mu_2^n} \frac{L_1(\mu^n_1)}{L_2(\mu_2^n)}
 > \bigg( \frac{\mu_1^{1+\ep }}{\mu_2^{1-\ep }}\bigg)^n\to  \infty
 \qquad \text{as }   n \to  \infty.
$$
Hence $f_n(0)<g_n(0)$ for all large enough $n$. Together with~\eqref{PT} that means
that  $F(s)<G(s)$ for all large enough $s$. Therefore $T_1$ will be ultimately dominant
by virtue of Theorem~\ref{t6}.
\end{proof}

Note that, in the case of critical BGW processes with offspring pgf's $f$ and $g$, a
result analogous to the first part of the above holds true: if the variance of $f$ is
lower than that of~$g$, then  $f_n(0)>g_n(0)$ for all large enough~$n$ (this follows,
say, from relation~(1) on p.~19 of~\cite{AN}).

The next results provides an upper bound for the number of dominance changes.

\begin {theorem}
 \label{t9}
Let $\{Z_1(n)\}$ and $\{Z_2(n)\}$ be two subcritical BGW processes with respective
offspring pgf's $f$ and $g$ and times to extinction $T_1$ and $T_2$. Then the number of
times the difference $f_n(0)-g_n(0),$ $n=1, 2,\ldots,$ changes sign and the number of
dominance changes for unmixed $n$-assemblies composed of copies of $T_1$ and $T_2$ are
no greater than the number of zeroes of $f(\theta)-g(\theta),$ $\theta \in [0,1)$.
\end{theorem}

\begin{proof}
Suppose that $f(\theta)\geq g(\theta)$ for $\theta\in [a,b) \subset [0,1]$ and
$f_{m}(0)\geq g_{m}(0)$ for some $m$ such that $f_{m}(0)\in[a,b)$. Then clearly
$f_{m+i+1}(0)\geq g_{m+i+1}(0)$ for all $i\geq 0$ such that $f_{m+i}(0) \in [a,b)$. Now
supposing that, for some $c \in (b, 1]$, $f(\theta)\leq g(\theta)$ for $\theta\in
[b,c)$ and $f_{n}(0)\leq g_{n}(0)$ for some $n$ such that $f_{n}(0)\in [b,c)$, we
likewise deduce that $f_{n+j+1}(0)\leq g_{n+j+1}(0)$ for all $j$ such that $f_{n+j}(0)
\in [b,c)$. It is easily seen that the above implies that, for any interval $[a,c)
\subset [0,1]$ with at most one zero of $f(\theta)-g(\theta)$, there will be at most
one change of sign in the finite subsequence  $\{f_{n}(0)-g_{n}(0): f_{n}(0)\in [a,c),
\ n\in \Nn\}$.

Suppose that $f(\theta)-g(\theta)$ has $\nu^*<\infty$ zeroes in $[0,1)$. Then,
subdividing  $[0,1)$ into subintervals with one zero in each, it follows that the total
number of changes of sign in the sequence $\{f_{n}(0)-g_{n}(0): \ n\in \Nn\}$ will not
exceed $\nu^*$.

Now let $F$ and $G$ be the respective distribution functions of $T_1$ and $T_2$, i.e.
$F(s)=f_{\lfloor s\rfloor}(0)$ and $G(s)=g_{\lfloor s\rfloor}(0)$. Then the number of
connected components of  $S^+:=\{s:F(s)>G(s)\}=\{s:f_{\lfloor s\rfloor}(0)>g_{\lfloor
s\rfloor}(0), s\in \mathbb{R_+}\}$ and $S^-:=\{s:F(s)<G(s)\}=\{s:f_{\lfloor
s\rfloor}(0)<g_{\lfloor s\rfloor}(0), s\in \mathbb{R_+}\}$ will not exceed $\nu^*$.
Applying Theorem~\ref{t7} completes the proof.
\end{proof}

Suppose now that the dynamics of the populations of the two species are modeled by
subcritical branching processes $\{Z_1(n)\}$ and $\{Z_2(n)\}$ with offspring pgf's $f$
and $g$, respectively, such that
\begin{equation}
 \label{cross}
f(\theta)-g(\theta)\quad \text{has at most one zero in\ } [0,1).
\end{equation}
This is the case, for instance, when both pgf's are quadratic:
$f(\theta)=p_0+p_1\theta+p_2\theta^2$ and $g(\theta)=q_0+q_1\theta+q_2\theta^2$ with
$p_0+p_1+p_2=q_0+q_1+q_2=1$.

The following result describes the behaviour of $M(k,n-k)$, $0\leq k\leq n$, depending
on agreement of the signs of the differences in respective offspring means and mean
times to extinction.

\begin{corollary}
 Assume that subcritical pgf's $f$ and $g$ satisfy condition~\eqref{cross}.
\begin{enumerate}
\item[{\rm(i)}]
If $(\mu_1-\mu_2)(\exn T_1-\exn T_2)>0$  then we have a dominant case, with $M(n,0)$
largest when $\mu_1>\mu_2$ and $M(0,n)$ largest otherwise.
\item[{\rm(ii)}]
If $(\mu_1-\mu_2)(\exn T_1-\exn T_2)<0$  then we have a non-dominant but ultimately
dominant case with exactly one sign change in the sequence $M(n,0)-M(0,n)$, $n=1, 2,
\ldots$
\end{enumerate}
\end{corollary}

\begin{proof}
(i)~Without loss of generality assume that $\mu_1>\mu_2$. Then, on the one hand, $T_1$
is ultimately dominant by Theorem~\ref{t8}, so that
$$
f_n(0)-g_n(0)<0 \qquad \text{for all large enough\ } n.
$$
On the other hand, due to condition~\eqref{cross}, by Theorem~\ref{t9} the sequence
$f_n(0)-g_n(0),$ $n=1, 2, \ldots,$ can change sign at most once. Therefore there exists
an $n_0\geq 0$ such that
\begin{align*}
 f_n(0)&-g_n(0)\geq 0, \qquad n\leq n_0,\\
 f_n(0)&-g_n(0)<0,     \qquad n>n_0,
\end{align*}
which clearly implies that $$\sum_{n\geq j}\big(f_n(0)-g_n(0)\big)<0, \qquad j>n_0,$$
whereas for $j\leq n_0$ one has
$$\sum_{n\geq j}\big(f_n(0)-g_n(0)\big)\leq \sum_{n\geq 0}\big(f_n(0)-g_n(0)\big)=\exn
T_2-\exn T_1<0$$ by our assumption. Hence $T_1$ is dominant by virtue of
Theorem~\ref{t4}.

\medskip

(ii)~We again assume that $\mu_1>\mu_2$, which implies that $T_1$ is ultimately
dominant owing to Theorem~\ref{t8}. However, $$\exn T_1=M(1,0)<M(0,1)=\exn T_2$$ by
assumption, so there is at least one change of sign in $M(n,0)-M(0,n)$. Finally, as
$f(\theta)-g(\theta), \ \theta \in [0,1)$, has at most one zero, Theorem~\ref{t9}
asserts that there cannot be more than one change of sign in the sequence.
\end{proof}

\subsubsection*{Acknowledgements.}
The first author would like to thank K. Hamza for many useful discussions. The research
was supported by the ARC Centre of Excellence for Mathematics and Statistics of Complex
Systems (MASCOS) and ARC grant DPO451657.

 \end{document}